\documentclass[12pt]{article}
\usepackage{amsfonts}
\usepackage{amssymb}

\usepackage{amsmath,url,epsfig,graphicx}
\newtheorem{lem}{Lemma}
\newtheorem{theorem}{Theorem}
\newtheorem{prop}{Proposition}

\title{Bi-quartic parametric polynomial minimal surfaces}

\author{Ognian Kassabov \and Krassimira Vlachkova}
\date{}
\begin{document}
\maketitle

\begin{abstract}
\noindent Minimal surfaces with isothermal parameters admitting B\'{e}zier representation were studied
by Cos\'{\i}n and Monterde. They showed that, up to an affine transformation,
the Enneper surface is the only bi-cubic isothermal minimal surface.
Here we study bi-quartic isothermal minimal surfaces and establish
the general form of their generating functions in the Weierstrass representation formula.
We apply an approach proposed by Ganchev to compute the normal curvature and
show that, in contrast to the bi-cubic case, there is a variety of
bi-quartic isothermal minimal surfaces. Based on the B\'{e}zier\  representation we establish some geometric properties of the bi-quartic harmonic surfaces. Numerical experiments are visualized and presented to illustrate and support our results.
\end{abstract}


\section{Introduction}

Minimal surfaces have recently become subject of intensive study in physical and biological sciences, e.g. materials science and molecular engineering which is due to its area minimizing property.  They are used in modeling physical phenomena as soap films, block copolymers, protein folding, solar cells, nanoporous membranes, etc. Minimal surfaces found applications also in architecture, CAGD, and computer graphics where B\'{e}zier\
polynomials and splines are widely used to efficiently describing,  representing and visualizing 3D objects. Hence it is important to know minimal surfaces in polynomial form of lower degrees. Bi-cubic
polynomial minimal surfaces are studied in \cite{C-M}.  Polynomial surfaces of degree 5
and 6 are studied in \cite{X-W-1} and \cite{X-W-2} where some interesting surfaces are described and their properties are examined. Examples of
polynomial minimal surfaces of arbitrary degree are presented in \cite{X-W-3}.

In section 3 we specify the result of Cos\'in and Monterde \cite{C-M} concerning
bi-cubic polynomial minimal surfaces. We note that their proof that these surfaces
coincide up to affine transformation with the classical Enneper surface concerns
the case of surfaces in {\it isothermal} parameters. In section 4 we consider
an analogous problem for polynomial surfaces defined by charts ${\bf x}(u,v)$ of
degree 4 on both $u$ and $v$. It turns out that they are
more various than the  bi-cubic ones, so they may be more useful in computer
graphics. In this paper we use a new approach to minimal surfaces proposed
by Ganchev \cite{GG} as well as the method from \cite{OK} to obtain a parametrization
of the surface in canonical parameters.

In section 5 we consider bi-quartic harmonic B\'{e}zier\  surfaces. We show that for a special choice of nine boundary control points
the corresponding harmonic B\'{e}zier\  surface is uniquely determined and is symmetric
with respect to one of the coordinate planes $Oxy$, $Oxz$, and $Oyz$. Based on
the B\'{e}zier representation we apply computer modeling and visualization
tools to illustrate and support our results.

\section{Preliminaries}\label{Preliminaries}
Let $S$ be a regular surface. Then $S$ is locally defined by a chart
$$
	{\bf x=x}(u,v)    \qquad (u,v) \in U \subset \mathbb R^2  \ .
$$
As usual we denote by $	{\bf x}_u$,	${\bf x}_v$, ${\bf x}_{uu}$,... the partial
derivatives of the vector function $ {\bf x}(u,v)$. Then the coefficients of
the first fundamental form are given by the inner products
$$
	E={\bf x}_u^2 \ , \qquad F={\bf x}_u{\bf x}_v \ , \qquad G={\bf x}_v^2
$$
and the unit normal  is
$$
	{\bf U} = \frac{{\bf x}_u\times {\bf x}_v}{|{\bf x}_u\times {\bf x}_v|}\cdot
$$
Then the coefficients of the second fundamental form are defined by
$$
	L={\bf U\,x}_{uu}, \quad 	M={\bf U\,x}_{uv}, \quad 	N={\bf U\,x}_{vv}.
$$
The Gauss curvature $K$ and the mean curvature $H$ of $S$ are given by
$$
  K=\frac{LN-M^2}{EG-F^2}, \quad	H=\frac{E N-2F M+G L}{2(EG-F^2)},
$$
respectively. Note that the Gauss curvature and the mean curvature
of a surface do not depend on the chart.
The surface $S$ is said to be {\it minimal} if its mean curvature
vanishes identically. In this case the Gauss curvature is negative and
the normal curvature of $S$ is the function $\nu=\sqrt{-K}$, see \cite{GG}.

We say that the chart ${\bf x}(u,v)$ is {\it isothermal} or that
the parameters $(u,v)$ are isothermal if $E=G$, $F=0$.
It is always possible to change the parameters $(u,v)$ so that the
resulting chart be isothermal. We note however that this change of the
parameters to isothermal ones is in general nonlinear.

When the chart is isothermal it is possible to use complex functions
to investigate it. We shall explain briefly this. Namely let $f(z)$ and
$g(z)$ be two holomorphic functions (actually sometimes they are taken
meromorphic). Define the Weierstrass complex curve $\Psi(z)$ by
\begin{equation}\label{2.1}
	\Psi(z)=\int_{z_0}^z \,
	\left( \frac12 f(z)(1-g^2(z)),\frac i2f(z)(1+g^2(z)),f(z)g(z) \right) \, dz\ .
\end{equation}
Then  $\Psi(z)$ is a minimal curve, i.e. $(\Psi'(z))^2=0$, and its real
and imaginary parts
$$
	{\bf x}(u,v)= Re\, \Psi(z) \quad {\rm and} \qquad {\bf y}(u,v)=Im\, \Psi(z)
$$
are minimal charts.
Moreover, they are isothermal and are harmonic functions
(i.e. $\Delta {\bf x}=0 $, $\Delta {\bf y}=0 $, where $\Delta$ is the Laplace
operator) as the real and complex part of a holomorphic function.
Conversely, every minimal surface can be defined at least locally in this way.
Of course a minimal surface can be generated by the Weierstrass formula with
different pairs of complex functions $f(z)$, $g(z)$.

It is easy to see that the coefficients of the first fundamental form
of a chart defined via the Weierstrass formula with functions $f(z)$, $g(z)$
are given by
\begin{equation}\label{2.2}
	E = G = \frac14 | f | ^2 (1 + | g | ^2)^2,  \quad F=0.
\end{equation}
The normal curvature is computed to be
\begin{equation}\label{2.3}
	\nu=\frac{4|g'|}{|f|(1+|g|^2)^2},
\end{equation}
see \cite{G-A-S}, Theorem 22.33.

Recently Ganchev \cite{GG} has proposed a new approach to minimal surfaces.
Briefly speaking he introduces special parameters called {\it canonical principal parameters}. A chart
with such parameters is isothermal. Moreover, the coefficients of the two fundamental forms
are given by
$$
\begin{array}{lll}
	E=\frac1\nu, &\ F=0, &\  G=\frac1\nu,\\[1ex]
	L=1, &\ M=0, &\ N=-1.
\end{array}
$$
His idea  leads to the fact that the real part of the minimal curve
\begin{equation}\label{gan1}
	\Phi(w)=-\int_{z_0}^z\left( \frac12\, \frac{1-\tilde g^2(w)}{\tilde g'(w)},\frac{i}2\,
	                             \frac{1+\tilde g^2(w)}{\tilde g'(w)},\frac{\tilde g(w)}{\tilde g'(w)} \right)\,dw
\end{equation}
is a minimal surface in canonical principal parameters. Note that
this is the Weierstrass formula with $f(z)=-1/\tilde g'(z)$, $g(z)=\tilde g(z)$.

We shall use also the following theorems:

{\bf Theorem A.} \cite{GG} If a surface is parametrized with canonical principal parameters,
then the normal curvature satisfies the equation
\begin{equation}\label{diffeq1}
	\Delta \ln \nu+2\nu=0.
\end{equation}
Conversely, for any solution $\nu(u,v)$ of equation \eqref{diffeq1} there exists
a {\bf unique} (up to position in the space) minimal surface with normal
curvature $\nu(u,v)$, where $(u,v)$ are canonical principal parameters.\hfill$\Box$

{\bf Theorem B.} \cite{OK} Let the minimal surface $S$ be defined by the
real part of \eqref{2.1}. Any solution of the differential equation
\begin{equation}\label{2.4}
	(z'(w))^2=-\frac1{f(z(w))g'(z(w))}
\end{equation}
defines a change of the isothermal parameters of $S$ to canonical principal
parameters. Moreover, the function $\tilde g(z)$ that defines $S$ via the Ganchev's
formula \eqref{gan1} is given by $\tilde g(w)=g(z(w))$.\hfill$\Box$

\vspace{0.2 cm}
The canonical principal parameters $(u,v)$ are determined
uniquely up to the changes
$$
	\begin{array}{l}
		u=\varepsilon\bar u+a,\\
		v=\varepsilon\bar v+b,
	\end{array} \quad \varepsilon=\pm 1, \ a=const.,\ b=const.
$$

\section{Bi-cubic minimal surfaces}

By investigating minimal B\'ezier surfaces Cos\'in and Monterde \cite{C-M} formulate
that any bi-cubic minimal surface defined by
$$
	{\bf x}(u,v)=\left(\sum_{i,j=0}^3 a_{ij}u^iv^j,\sum_{i,j=0}^3 b_{ij}u^iv^j,\sum_{i,j=0}^3 c_{ij}u^iv^j\right)
$$
is, up to affine transformation in the space, actually an affine reparametrization
of the classical Enneper surface
$$
	{\bf enneper}(u,v)=\frac12\left( u-\frac{u^3}{3}+uv^2,-v+\frac{v^3}{3}-u^2v,u^2-v^2 \right).
$$
(This chart of the Enneper surface is obtained from the Weierstrass formula with $f(z)=1$, $g(z)=z$.)
We note that their proof actually refers to bi-cubic {\bf isothermal} minimal charts.
Indeed, as we mentioned in Section \ref{Preliminaries}, the change of parameters to
isothermal ones is in general nonlinear. Below we give a simple example of a
bi-cubic minimal chart that can not be transformed by an affine transformation into
an isothermal one.

\vspace{0.2cm}
{\bf Example.} Consider the bi-cubic chart
\begin{equation}\label{3.1}
	{\bf x}(u,v)=\frac12 \left( u v - \frac{u^3 v^3}3 + u v^3, -v + \frac{v^3}3 - u^2 v^3,  u^2 v^2-v^2 \right).  \end{equation}
It can be shown by direct computation that this chart defines a minimal surface. Of course
we can simply remark that this is a reparametrization of ${\bf enneper}(u,v)$ with
$u$ replaced by $uv$, so the mean curvature vanishes identically.

Let us make an affine transformation of the parameters $(u,v)$:
$$
	\begin{array}{l}
		u=a_1 \bar u+b_1\bar v+c_1 \\
		v=a_2 \bar u+b_2\bar v+c_2
	\end{array}
$$
with nonzero Jacobian, i.e.
\begin{equation}\label{3.2}
	J=a_1b_2-a_2b_1 \ne 0.
\end{equation}
We shall try to determine the coefficients $a_i, b_i, c_i$, so that
the chart
$$
	\bar{\bf x}(\bar u,\bar v)={\bf x}(a_1 \bar u+b_1\bar v+c_1,a_2 \bar u+b_2\bar v+c_2)
$$
be isothermal. Actually we shall see what follows only from $\overline F=0$.
A direct computation shows that $\overline F=\bar{\bf x}_{\bar u}\cdot\bar{\bf x}_{
\bar v}=\frac14F_1\cdot F_2$, where
$$
\begin{array}{ll}
	F_1=& \Bigl(1 + (c_2 + a_2  {\bar u} + b_2  {\bar v})^2 \bigl(1 + (c_1 + a_1  {\bar u} + b_1  {\bar v})^2\bigr)\Bigr)^2,\\
F_2=& a_2^2 b_1 {\bar u} (2 a_1 {\bar u} + b_1 {\bar v}+c_1) +
  a_1 ( b_2 {\bar v}+c_2) \bigl(b_1 c_2 + b_2 (a_1 {\bar u} + 2 b_1 {\bar v}+c_1)\bigr)\\
  & +a_2 \Bigl(b_1 c_2 (3 a_1 {\bar u} + b_1 {\bar v}+c_1) + b_2 \bigl(1 + c_1^2 + 3 a_1 c_1 {\bar u} + 2 a_1^2 {\bar u}^2 + 3 b_1 (c_1 + 2 a_1 {\bar u}) {\bar v} +
        2 b_1^2 {\bar v}^2\bigr)\Bigr).
     \end{array}
$$
Since $F_1$ is positive, the vanishing of $\overline F$ implies $F_2=0$. Hence the
coefficients in  $F_2$ must  be zero. In particular the coefficients of
${\bar u}^2$, ${\bar v}^2$, ${\bar u}{\bar v}$ are
$$
\begin{array}{l}
	 a_1 a_2(a_2 b_1 +  a_1 b_2)=0,\\
	 b_1b_2(a_2 b_1 +  a_1 b_2)=0,\\
	(a_2 b_1 + a_1b_2)^2+4 a_1 a_2 b_1 b_2 =0.
\end{array}
$$
These equations imply immediately
\begin{equation}\label{3.3}
  a_1a_2b_1b_2=0, \quad a_2 b_1 + a_1b_2=0.
\end{equation}
Let e.g. $a_1=0$. From \eqref{3.3} it follows $a_2b_1=0$, which contradicts \eqref{3.2}.
So it is impossible to make an affine transformation of the parameters in
\eqref{3.1} to obtain an isothermal chart.

{\bf Remark.} {In view of the above notes the problem of existing bi-cubic
minimal surfaces different from the Enneper one is still open. More
generally it will be interesting to obtain a method for finding
polynomial minimal non-isothermal charts.}

\section{Bi-quartic minimal surfaces in isothermal parameters}

In this section we examine minimal surfaces represented by
isothermal polynomial charts of degree 4 in both
$u,v$. We may expect that there exists more than one such surface, but it
is interesting to know ``how many'' are there.

So consider the chart
$$
	{\bf x}(u,v)=\sum_{i,j=0}^4 {\bf v}_{ij}  u^iv^j,
$$
where ${\bf v}_{ij} =(a_{ij},b_{ij},c_{ij})$ are vectors in $\mathbb R^3$.
Using $F=0$ and looking on its coefficient of $u^7v^7$ we obtain
${\bf v}_{44}=0$. Analogously we derive consecutively ${\bf v}_{43}=0$,
${\bf v}_{34}=0$, ${\bf v}_{42}=0$, ${\bf v}_{24}=0$, ${\bf v}_{41}=0$, ${\bf v}_{14}=0$,
${\bf v}_{33}=0$, ${\bf v}_{32}=0$, ${\bf v}_{23}=0$.

It is known that any minimal isothermal chart is harmonic, see e.g. \cite{G-A-S}.
In our case this implies
$$
\begin{array}{lll}
	{\bf v}_{02}=-{\bf v}_{20}, &\ {\bf v}_{21}=-3{\bf v}_{03},&\ {\bf v}_{22}=-6{\bf v}_{40},\\
	{\bf v}_{12}=-3{\bf v}_{30}, &\ {\bf v}_{13}=-{\bf v}_{31},&\ {\bf v}_{04}={\bf v}_{40}.
\end{array}
$$
Substituting these in $F$ and looking on the coefficients of $u^6$
and $u^5v$  we obtain also
${\bf v}_{31}{\bf v}_{40}=0$, ${\bf v}_{31}^2-16{\bf v}_{40}^2=0$.
If ${\bf v}_{40}=\bf o$ the chart is not of degree 4. So we assume
${\bf v}_{40}\ne \bf o$. Up to position in space and symmetry we may take
$$
	{\bf v}_{40} = (p, 0, 0),\quad {\bf v}_{31} = (0, 4p, 0),\quad {\rm where} \ p\ne 0.
$$
Now the coefficients of $u^5$  and $u^4v$ in $F$ give
$b_{03} + a_{30} = 0$ and $a_{03} - b_{30} = 0$.
Using this we can calculate the derivatives
${\bf x}_u$ and ${\bf x}_v$.
Let the functions $f(z)$ and $g(z)$ give the Weierstrass representation
of the surface. Denote by $(\phi_1,\phi_2,\phi_3)$ the derivative of $\Psi$.
Then
$
	\Psi'=(\phi_1,\phi_2,\phi_3)={\bf x}_u-i {\bf x}_v.
$
In our case a direct computation shows
$$
	\begin{array}{l}
  		\phi_1=	-i a_{01} +  a_{10} + (u + i v) \Bigl(-i a_{11} +
            2 a_{20} + (u + i v) \bigl(3 a_{30} + 3 i b_{30} + 4 p( u +  i  v)\bigr)\Bigr),\\
   		\phi_2= -i b_{01} + b_{10} -i ( u +i v) \Bigl(b_{11} +
          2 i b_{20} + (u + i v) \bigl(3 a_{30} + 3 i b_{30} + 4 p (u +  i v)\bigr)\Bigr),\\
      \phi_3= -i c_{01} + c_{10} + \bigl(-i c_{11} + 2 c_{20} + 3 (i c_{03} + c_{30}) (u + i v)\bigr) (u + i v).
	\end{array}
$$
On the other hand the Weierstrass formula implies easily
$$
	f(z)=\phi_1-i\phi_2,\quad g(z)=\frac{\phi_3}{\phi_1-i\phi_2}.
$$
Hence we derive
\begin{eqnarray*}
	f(z)&=&-i a_{01} - b_{01} + a_{10} - i b_{10} + (-i a_{11} - b_{11} + 2 a_{20} - 2 i b_{20}) z,\\[4pt]
	g(z)&=&\frac{ c_{01} + i c_{10} + (c_{11} + 2 i c_{20} + 3 (- c_{03} + i c_{30}) z) z}
	{a_{01} -  i b_{01} + i a_{10} + b_{10} + (a_{11} - i b_{11} + 2 i a_{20} + 2 b_{20}) z}\cdot
\end{eqnarray*}
Consequently we have obtained that for some complex constants $A$ and $B$
$$
	f(z)=Az+B, \qquad   g(z)=\frac{P_2(z)}{Az+B},
$$
where $P_2(z)$ is a polynomial of degree at most 2. Suppose $A=0$, i.e. $f(z)$ is a
constant. Then the derivative
$$
		(\phi_1,\phi_2,\phi_3)=\left( \frac12 f(z)(1-g^2(z)),\frac i2 f(z)(1+g^2(z),f(z)g(z) \right)
$$
is of degree 2 or 4, so the chart ${\bf x}(u,v)$ is of degree 3 or 5, which is not our case. So $A\ne 0$.
Since $
	\phi_1=\frac12f(z)(1-g^2(z))
$
is a polynomial then $Az+B$ divides
$P_2(z)$. Hence $g(z)=Cz+D$,  where $C\ne 0$. We have proved the following
\begin{theorem}\label{th1}
{Any bi-quartic parametric polynomial minimal surface in isothermal
parameters is generated by the Weierstrass formula with the functions
$$f(z)=Az+B,\ g(z)=Cz+D,\ \mbox{where}\ A\not=0,\ C\not=0.$$}
\end{theorem}
Further, we are interested which of the functions in Theorem \ref{th1} generate different surfaces.
Denote by ${\bf x_0}(u,v)$ the chart defined as the real part of the Weierstrass
minimal curve with functions
$f(z)=z$, $g(z)=z$ and the corresponding surface by $S_0$.
Denote also by ${\bf x}(u,v)$ the chart with
generating functions $f(z)=A z$, $g(z)=C z$ for arbitrary
nonzero complex numbers $A$, $C$, and the corresponding surface by $S$.
Using \eqref{2.2} and \eqref{2.3} we can see that the nonzero coefficients
of the first fundamental form and the normal curvature of
${\bf x_0}(u,v)$ are respectively
$$
	E_0=G_0=\frac14 (u^2 + v^2) (1 + u^2 + v^2)^2,\ \mbox{and\ hence}\ 
\nu_0=\frac4{\sqrt{u^2 + v^2} (1 + u^2 + v^2)^2}\cdot
$$
Obviously ${\bf x_0}(u,v)$ is not in canonical principal parameters.
We want to change the parameters $(u,v)$ to canonical principal ones.
Equation \eqref{2.4} has a solution
$$
	z= (3/2)^{2/3} (i\, w)^{2/3}.
$$
So according to Theorem B we change the variable $z$ by $ (3/2)^{2/3} (i\, z)^{2/3}$.
Now the  functions
$$
	\tilde f(z)=i\,\left(\frac32\right)^{1/3} (i\, z)^{1/3}, \quad \tilde g(z)=\left(\frac32\right)^{2/3} (i\, z)^{2/3}
$$
generate a chart ${\bf \tilde x_0}(u,v)$ in canonical principal
parameters and
$$
	\tilde\nu_0=\frac{4 \left(\frac23\right)^{2/3}}{(u^2 + v^2)^{1/3} \left(1 + \left(\frac32\right)^{4/3} (u^2 + v^2)^{2/3}\right)^2}\cdot
$$

\begin{figure}[hbtp]
\begin{minipage}[b]{2.5in}
\centering
\includegraphics[width=0.8\textwidth]{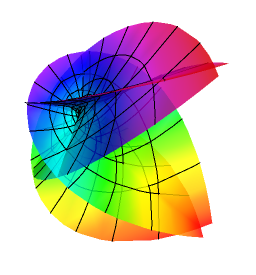}
  \centerline{\footnotesize Generating\ functions\ f(z)=10z,\ g(z)=z}
  \end{minipage}
~~~~~~~~~
\begin{minipage}[b]{2.5in}
\centering
\includegraphics[width=0.8\textwidth]{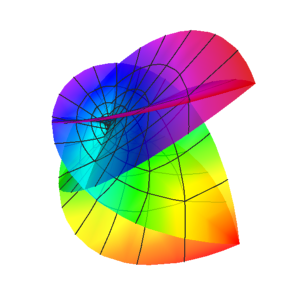}
    \centerline{\footnotesize Generating\ functions\ f(z)=z,\ g(z)=10z}
    \end{minipage}
\caption{Surfaces of type $S$ for different values of $A$ and $C$}\label{generating}
\end{figure}

Analogously ${\bf x}(u,v)$ is not in canonical principal parameters.
Using again Theorem B and changing the complex variable $z$ by
$$
	\left(\frac32\right)^{2/3} \left(\frac{i\, z}{\sqrt A \sqrt C}\right)^{2/3}
$$
we obtain a corresponding chart ${\bf \tilde x}(u,v)$ in canonical principal parameters.
According to \eqref{2.3} its normal curvature is
$$
	\tilde\nu = \frac{4 \left(\frac23\right)^{2/3} \left(\frac{|C|^2}{|A|}\right)^{2/3}}
	         {  (u^2 + v^2)^{1/3} \left(1 +  \left(\frac32\right)^{4/3} \left(\frac{| C |^2}{|A|}\right)^{2/3} (u^2 + v^2)^{2/3}\right)^2}\cdot
$$
The last formula implies that this is also the normal curvature
(in canonical principal parameters) of the surface $S_1$
generating via the Weierstrass formula by the functions
$$
	f_1(z)=\frac{|A|}{|C|^2}\,z, \quad g_1(z)=z.
$$
According to Theorem A the surfaces $S$ and $S_1$ coincide up to position in the
space. On the other hand, the Weierstrass formula implies that the surfaces $S_0$
and $S_1$ are homothetic. So for any nonzero complex numbers $A$, $B$ the surface
$S$ is, up to position in the space, homothetic to $S_0$.
Surfaces of type $S$ for different values of $A$ and $C$ are shown in Figure \ref{generating}.

\begin{figure}[hbtp]
\begin{minipage}[b]{1.8in}
\centering
\includegraphics[width=0.9\textwidth]{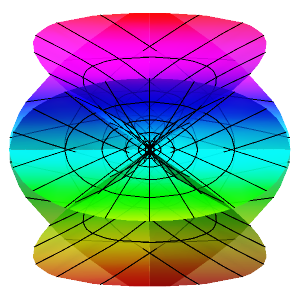}
  \centerline{\footnotesize $f(z)=z,\ g(z)=z$}
  \end{minipage}
~~~~~~~~~
\begin{minipage}[b]{1.8in}
\centering
\includegraphics[width=0.9\textwidth]{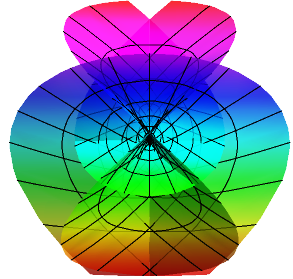}
    \centerline{\footnotesize $f(z)=z,\ g(z)=z+1/2$}
    \end{minipage}
    ~~~~~~~~~
\begin{minipage}[b]{1.8in}
\centering
\includegraphics[width=0.9\textwidth]{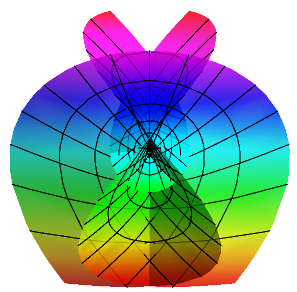}
    \centerline{\footnotesize $f(z)=z,\ g(z)=z+1$}
    \end{minipage}
    \caption{Comparison of bi-quartic minimal surfaces for generating functions $f(z)=z$ and $g(z)=z+a+i\,b$}\label{comparison}
    \end{figure}

Now we are interested whether the functions $f(z)=z$, $g(z)=z+a+i\,b$, $a,\,b\in\mathbb R$,
define a surface which is really different from $S_0$. The chart ${\bf x}(u,v)$ generated by these
functions is not in canonical principal parameters so according to Theorem B we
change the complex variable $z$ by $ (3/2)^{2/3} (i\, z)^{2/3}$. Then the functions
$$
	\tilde f(z)=i\, \left(\frac32\right)^{1/3} (i\, z)^{1/3}, \quad \tilde g(z)= \left(\frac32\right)^{2/3}(i\, z)^{2/3}+a + i\,b
$$
define a chart ${\bf \tilde x}(u,v)$ in canonical principal
parameters. Its normal curvature is
$$
	\tilde\nu=\frac{4\left(\frac23\right)^{2/3}}{ \sqrt{u^2 + v^2}
	( 1 + B(u,v)\overline B(u,v))^2},
$$
where $B(u,v)$ is
$$
  B(u,v)=a + i\,b + \left(\frac32\right)^{2/3} (i\,u -  v)^{2/3}.
$$
Comparing these functions for different values of $(a,b)$    we can say that
the resulting surfaces are different. In Figure \ref{comparison} are shown parts of the surface $S_0$,
obtained for $(a,b)=(0,0)$ (left), the surface obtained for $(a,b)=(0.5,0)$ (center), and the surface obtained for $(a,b)=(1,0)$ (right).

\section{Bi-quartic harmonic B\'{e}zier\  surfaces}

We consider bi-quartic tensor product B\'{e}zier\  surface defined by
\begin{equation}\label{eq1}
{\bf x}(u,v)=\sum_{i=0}^4\sum_{j=0}^4{\bf b}_{ij}B_i^4(u)B_j^4(v),
 \end{equation}
 where ${\bf b}_{ij},\ i,j=0,\ldots, 4 $ are the control points of ${\bf x}(u,v)$, and $B_i^4(u)$ are the Bernstein polynomials of degree $4$ defined for $0\leq u\leq 1$ by
$$B_i^4(u):=\binom{4}{i} u^i(1-u)^{4-i}, \qquad
\binom{4}{i}= \begin{cases}
                 \frac{4!}{i!(4-i)!}, & \text{for }\ i=0,\dots ,4, \\
                 0,                   & \text{otherwise.}
              \end{cases} $$
Recall that if ${\bf x}(u,v)$ is in isothermal parameters then ${\bf x}(u,v)$ is a minimal surface if and only if ${\bf x}(u,v)$ is a harmonic surface, i.e. $\Delta {\bf x}=0$. For a harmonic B\'{e}zier\  surface Monterde \cite{M} has proved that if we know the control points on two opposite boundaries except one corner point, e.g. nine points $\{{\bf b}_{0j}\}_{j=0}^4$ and $\{{\bf b}_{i4}\}_{i=0}^3$, then the remaining sixteen control points are fully determined. The proof\footnote{Monterde's proof is made for a surface of degree $(n,n)$, where $n\in\mathbb{N}$ is even. Here we consider the case $n=4$.} is based on the harmonic condition $\Delta {\bf x}=0$ which leads to a linear system that has a unique solution.
\begin{figure}[hbtp]
\centering
$\begin{array}{ccccc}
  {\bf b}_{40} &\circ &\circ &\circ &{\bf b}_{44}\\[1ex]
  {\bf b}_{30}&\circ &\circ &\circ &{\bf b}_{34}\\[1ex]
  {\bf b}_{20}&\circ &\circ &\circ &\circ\\[1ex]
{\bf b}_{10}&\circ &\circ&\circ &{\bf b}_{14}\\[1ex]
{\bf b}_{00}&\circ &\circ &\circ &{\bf b}_{04}
\end{array} $
  \caption{\small Input control points that fully determined bi-quartic harmonic B\'{e}zier\  surface}\label{f2}
\end{figure}
Here we assume that nine control points $\{{\bf b}_{i0}\}_{i=0}^4$, $\{{\bf b}_{i4}\}_{i=0;i\not=2}^4$ as shown on Figure \ref{f2} are given. Note that they differ from those used in \cite{M}. In Lemma \ref{lemma1} below we give expressions for the remaining control points through the given points. Then in Proposition \ref{pr1} we prove that in the case where the given points are symmetric with respect to any of the coordinate planes, then the corresponding harmonic B\'{e}zier\  surface is symmetric with respect to the same coordinate plane.
\begin{lem}\label{lemma1}
{Let nine control points ${\bf b}_{i0},\ i=0,\dots,4$; ${\bf b}_{i4},\ i=0,1,3,4$ be given. A bi-quartic B\'{e}zier\  surface \eqref{eq1} is harmonic if and only if the remaining sixteen control points satisfy}
\end{lem}
\begin{equation}\label{eq3}
\begin{array}{rcl}
{\bf b}_{01}&=&(25{\bf b}_{00}+8{\bf b}_{04}-40{\bf b}_{10}-8{\bf b}_{14}+36{\bf b}_{20}-16{\bf b}_{30}
+4{\bf b}_{34}+4{\bf b}_{40}-{\bf b}_{44}  )/12\\
{\bf b}_{41}&=&(4{\bf b}_{00}-{\bf b}_{04}-16{\bf b}_{10}+4{\bf b}_{14}+36{\bf b}_{20}-40{\bf b}_{30}
-8{\bf b}_{34}+25{\bf b}_{40}+8{\bf b}_{44}  )/12\\
{\bf b}_{11}&=&(17{\bf b}_{00}+7{\bf b}_{04}-14{\bf b}_{10}-4{\bf b}_{14}+18{\bf b}_{20}-8{\bf b}_{30}
+2{\bf b}_{34}+5{\bf b}_{40}+{\bf b}_{44}  )/24\\
{\bf b}_{31}&=&(5{\bf b}_{00}+{\bf b}_{04}-8{\bf b}_{10}+2{\bf b}_{14}+18{\bf b}_{20}-14{\bf b}_{30}
-4{\bf b}_{34}+17{\bf b}_{40}+7{\bf b}_{44}  )/24\\
{\bf b}_{02}&=&
(13{\bf b}_{00}+8{\bf b}_{04}-28{\bf b}_{10}-8{\bf b}_{14}+30{\bf b}_{20}-16{\bf b}_{30}+4{\bf b}_{34}+4{\bf b}_{40}-{\bf b}_{44}  )/6\\
{\bf b}_{42}&=&(4{\bf b}_{00}-{\bf b}_{04}-16{\bf b}_{10}+4{\bf b}_{14}
+30{\bf b}_{20}-28{\bf b}_{30}-8{\bf b}_{34}+13{\bf b}_{40}+8{\bf b}_{44}  )/6\\
{\bf b}_{12}&=&(11{\bf b}_{00}+7{\bf b}_{04}-20{\bf b}_{10}-4{\bf b}_{14}
+24{\bf b}_{20}-14{\bf b}_{30}+2{\bf b}_{34}+5{\bf b}_{40}+{\bf b}_{44}  )/12\\
{\bf b}_{32}&=&(5{\bf b}_{00}+{\bf b}_{04}-14{\bf b}_{10}+2{\bf b}_{14}
+24{\bf b}_{20}-20{\bf b}_{30}-4{\bf b}_{34}+11{\bf b}_{40}+7{\bf b}_{44}  )/12\\
{\bf b}_{03}&=&(14{\bf b}_{00}+19{\bf b}_{04}-32{\bf b}_{10}-16{\bf b}_{14}
+36{\bf b}_{20}-20{\bf b}_{30}+8{\bf b}_{34}+5{\bf b}_{40}-2{\bf b}_{44}  )/12\\
{\bf b}_{43}&=&(5{\bf b}_{00}-2{\bf b}_{04}-20{\bf b}_{10}+8{\bf b}_{14}
+36{\bf b}_{20}-32{\bf b}_{30}-16{\bf b}_{34}+14{\bf b}_{40}+19{\bf b}_{44}  )/12\\
{\bf b}_{13}&=&(5{\bf b}_{00}+7{\bf b}_{04}-8{\bf b}_{10}-{\bf b}_{14}
+9{\bf b}_{20}-5{\bf b}_{30}+2{\bf b}_{34}+2{\bf b}_{40}+{\bf b}_{44}  )/12\\
{\bf b}_{33}&=&(2{\bf b}_{00}+{\bf b}_{04}-5{\bf b}_{10}-2{\bf b}_{14}
+9{\bf b}_{20}-8{\bf b}_{30}-{\bf b}_{34}+5{\bf b}_{40}+7{\bf b}_{44}  )/12\\
{\bf b}_{21}&=&(3{\bf b}_{00}+{\bf b}_{04}-4{\bf b}_{10}+8{\bf b}_{20}-4{\bf b}_{30}+3{\bf b}_{40}
+{\bf b}_{44})/8\\
{\bf b}_{22}&=&(7{\bf b}_{00}+3{\bf b}_{04}-16{\bf b}_{10}+24{\bf b}_{20}-16{\bf b}_{30}+7{\bf b}_{40}
+3{\bf b}_{44})/12\\
{\bf b}_{23}&=&(7{\bf b}_{00}+5{\bf b}_{04}-16{\bf b}_{10}+4{\bf b}_{14}
+24{\bf b}_{20}-16{\bf b}_{30}+4{\bf b}_{34}+7{\bf b}_{40}+5{\bf b}_{44}  )/24\\
{\bf b}_{24}&=&({\bf b}_{00}-{\bf b}_{04}-4{\bf b}_{10}+4{\bf b}_{14}
+6{\bf b}_{20}-4{\bf b}_{30}+4{\bf b}_{34}+{\bf b}_{40}-{\bf b}_{44}  )/6.
\end{array}
\end{equation}

{\bf Proof.} It follows straightforward using the corresponding linear system from \cite{M}.
\hfill$\Box$

\begin{prop}\label{pr1}
{Let the given points ${\bf b}_{i0},\ i=0,\dots,4$; ${\bf b}_{i4},\ i=0,1,3,4$ be symmetric with respect to some of the  coordinate planes $Oxy$, $Oxz$, or $Oyz$.  Then the corresponding harmonic B\'{e}zier\  surface defined by Lemma \ref{lemma1} is symmetric with respect to the same plane. }
\end{prop}
{\bf Proof.} Let ${\bf b}_{ij}={\bf b}_{ij}(x_{ij},y_{ij},z_{ij})$ and assume that the given control points are symmetric with respect to the plane $Oxy$, i.e. ${\bf b}_{0k}$ and ${\bf b}_{1k}$ are symmetric points to ${\bf b}_{4k}$ and ${\bf b}_{3k}$, respectively, $k=0,4$, and ${\bf b}_{20}$ lies on $Oxy$. Then we have
\begin{equation}\label{eq2}
x_{0k}=x_{4k},\ x_{1k}=x_{3k},\ y_{0k}=y_{4k},\ y_{1k}=y_{3k},\ z_{0k}=-z_{4k},\ z_{1k}=-z_{3k}
\end{equation}
 for $k=0,4$, and $z_{20}=0$. To show that the harmonic B\'{e}zier\  surface defined by Lemma \ref{lemma1} is symmetric with respect to $Oxy$ it suffices to establish that its control points are symmetric with respect to $Oxy$. We need to establish that \eqref{eq2} holds for $k=1,2,3$ and $z_{2j}=0$ for $j=1,\dots ,4$. Next we verify that $x_{01}=x_{41}$, $y_{01}=y_{41}$, $z_{01}=-z_{41}$, and $z_{21}=0$. The analogous relations for the remaining control points follow in a similar way.

 From \eqref{eq3} and \eqref{eq2} we have
 $$\begin{array}{ll}
 x_{01}&=(25x_{00}+8x_{04}-40x_{10}-8x_{14}+36x_{20}-16x_{30}
+4x_{34}+4x_{40}-x_{44}  )/12\\
&=(29x_{00}+7x_{04}-56x_{10}-4x_{14}+36x_{20})/12=x_{41}.
\end{array}$$
 Analogous relation holds for $y_{01}$ and $y_{41}$. For the third coordinates $z_{01}$ and $z_{41}$ we obtain
$$z_{01}=(21z_{00}+9z_{04}-24z_{10}-12z_{14})/12=-z_{41}.$$
It remains to show that $z_{21}=0$. We have
\begin{eqnarray*}
z_{21}&=&(3z_{00}+z_{04}-4z_{10}+8z_{20}-4z_{30}+3z_{40}
+z_{44})/8\\
&=&\bigl(3(z_{00}+z_{40})+(z_{04}+z_{44})-4(z_{10}+z_{30})+8z_{20}\bigr)/8=0.
\end{eqnarray*}
The case where the nine given points are symmetric with respect to the other coordinate planes is treated analogously.
\hfill$\Box$

\begin{figure}[hbtp]
\begin{minipage}[b]{2.5in}
\centering
\includegraphics[width=.95\textwidth]{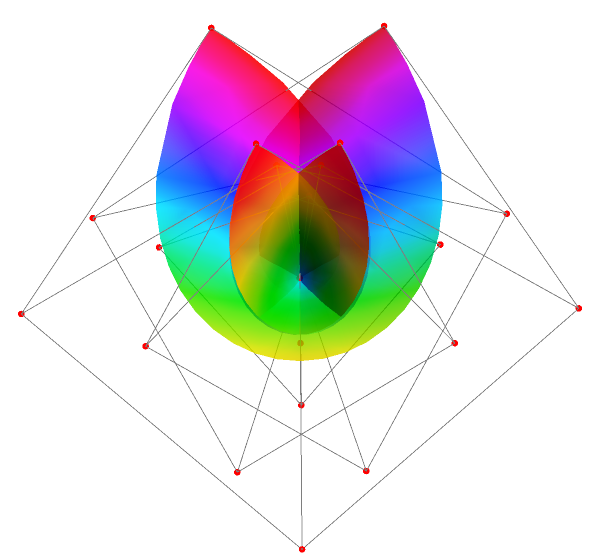}
                  \end{minipage}
~~~~~~~~~
\begin{minipage}[b]{2.5in}
\centering
\includegraphics[width=.8\textwidth]{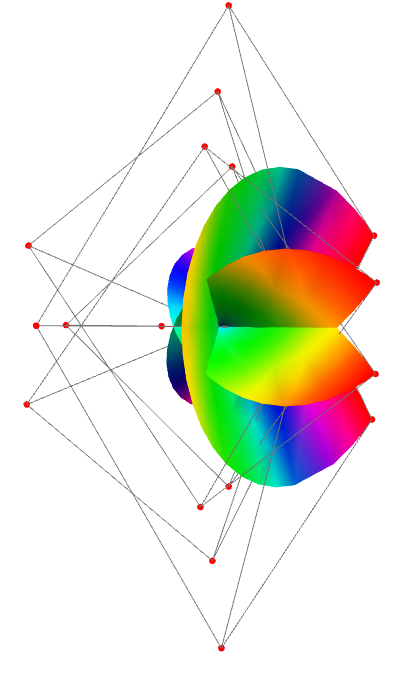}
                  \end{minipage}
                  \caption{\small Symmetry of a harmonic bi-quartic B\'{e}zier\  surface with respect to $Oxz$. The surface is shown from two different viewpoints.}
  \label{Bezier1}
\end{figure}

A bi-quartic harmonic B\'{e}zier\  surface which is symmetric with respect to $Oxz$ is shown from two different viewpoints in Figure \ref{Bezier1}. Its control points are presented in Table \ref{tab1}. We note that they are obtained from the minimal bi-quartic B\'{e}zier\  surface with generating functions $f(z)=z$, $g(z)=z-1$. Hence, the surface in Figure \ref{Bezier1} is harmonic minimal B\'{e}zier\  surface.
\arraycolsep=9pt
\arraycolsep=4pt
\begin{table}
\begin{tabular}{ccccc}
\hline\\
$(\frac{512}{3},-\frac{176}{3},\frac{128}{3})$& $(\frac{128}{3}, \frac{424}{3}, \frac{152}{3})$&$ (-128, \frac{128}{3}, \frac{32}{3})$&$ (-\frac{128}{3}, -\frac{296}{3}, -\frac{104}{3})$& $(\frac{256}{3}, -\frac{80}{3},-\frac{128}{3})$\\[1ex]
$(-\frac{64}{3}, -\frac{536}{3}, -\frac{88}{3})$&$ (\frac{224}{3}, -\frac{44}{3}, \frac{32}{3})$&$ (0, 0, \frac{8}{3})$&$ (\frac{160}{3}, -\frac{20}{3}, -\frac{32}{3})$& $(\frac{64}{3}, \frac{280}{3}, \frac{40}{3})$\\[1ex]
$(-\frac{512}{3}, 0, -\frac{160}{3})$&$ (0, 0, -\frac{8}{3})$&$ (-\frac{128}{3}, 0, 0)$&$ (0, 0, -\frac{8}{3})$& $ (-\frac{256}{3}, 0, 32)$\\[1ex]
$(-\frac{64}{3}, \frac{536}{3}, -\frac{88}{3})$&$ (\frac{224}{3}, \frac{44}{3}, \frac{32}{3})$&$ (0, 0, \frac{8}{3})$&$ (\frac{160}{3}, \frac{20}{3}, -\frac{32}{3})$&$ (\frac{64}{3}, -\frac{280}{3}, \frac{40}{3})$\\[1ex]
$(\frac{512}{3}, \frac{176}{3}, \frac{128}{3})$&$ (\frac{128}{3}, -\frac{424}{3}, \frac{152}{3})$&$ (-128, -\frac{128}{3}, \frac{32}{3})$&$ (-\frac{128}{3}, \frac{296}{3}, -\frac{104}{3})$&$ (\frac{256}{3}, \frac{80}{3}, -\frac{128}{3})$\\[4pt]
\hline
\end{tabular}
\caption{Control points of a harmonic bi-quartic B\'{e}zier\  surface that are symmetric with respect to $Oxz$}
\label{tab1}
\end{table}

\section{Conclusions and Future Work}

In this paper we characterize all bi-quartic parametric polynomial minimal surfaces by their generating functions using the Weierstrass formula. We also consider the bi-quartic harmonic B\'{e}zier\ surfaces and establish their symmetry with respect to any of the coordinate planes. We present numerical experiments
and give examples. A possible direction for future work is to extend our results for minimal surfaces
of higher degrees.

\vspace{.2cm}
{\bf Acknowledgments.}  This work was partially supported by the Bulgarian National Science Fund
under Grant No. DFNI-T01/0001.

\bigskip\noindent
{\sc Ognian Kassabov} \smallskip\\
Department of Mathematics and Informatics\\
University of Transport ``Todor Kableshkov''\\
158 G. Milev Str. \\
1574 Sofia, BULGARIA\\
email: \texttt{okassabov@abv.bg}\\

\noindent
{\sc Krassimira Vlachkova} \smallskip\\
Faculty of Mathematics and Informatics\\
Sofia University\\
5 James Bourchier Blvd. \\
1164 Sofia, BULGARIA\\
email: \texttt{krassivl@fmi.uni-sofia.bg}\\
\end {document}